\documentclass{amsart}
\usepackage{amsfonts}
\usepackage{amssymb}
\usepackage{xypic}

\newcommand{\Z}{{\mathbb Z}}
\newcommand{\C}{{\mathbb C}}
\renewcommand{\P}{{\mathbb P}}
\newcommand{\QH}{\text{\sl QH}}

\newcommand{\Gr}{\operatorname{Gr}}
\newcommand{\Fl}{\operatorname{F\ell}}

\newcommand{\bull}{{\sssize \bullet}}

\newtheorem{lemma}{Lemma} 
\newtheorem{prop}{Proposition} 
\newtheorem{thm}{Theorem}

\newcommand{\reflemma}[1]{Lemma~\ref{#1}}

\newcommand{\refprop}[1]{Proposition~\ref{#1}}
\newcommand{\refsec}[1]{Section~\ref{#1}}

\begin{document}

\title{Direct proof of the quantum Monk's formula}
\author{Anders Skovsted Buch}
\address{Massachusetts Institute of Technology \\
  Building 2, Room 275 \\
  77 Massachusetts Avenue \\
  Cambridge, MA 02139
}
\date{\today}
\email{abuch@math.mit.edu}
\thanks{The author was partially supported by NSF Grant DMS-0070479}
\maketitle

\section{Introduction}

The quantum Monk's formula of Fomin, Gelfand, and Postnikov
\cite{fomin.gelfand.ea:quantum} gives an explicit rule for multiplying
by a codimension one Schubert class in the (small) quantum cohomology
ring of a flag variety $SL_n/B$.  In the present paper we give a
direct geometric proof of this formula which relies only on classical
Schubert calculus.  In particular, no compactifications of moduli
spaces are required.  Our proof uses an adaption of the ideas from
\cite{buch:quantum} where we give a similar proof of the quantum Pieri
formula for Grassmann varieties.

Since the quantum cohomology ring of a flag variety is generated by
the codimension one Schubert classes, the quantum Monk's formula
uniquely determines this ring as well as the associated Gromov-Witten
invariants.  Thus, if associativity of quantum cohomology is granted
\cite{ruan.tian:mathematical, kontsevich.manin:gromov-witten,
  fulton.pandharipande:notes}, we obtain a complete elementary
understanding of this ring.

The presentation due to Givental, Kim, and Ciocan-Fontanine
\cite{givental.kim:quantum, kim:quantum*2, ciocan-fontanine:quantum*1}
and Ciocan-Fontanine's formula for the quantum Schubert classes given
by certain cyclic permutations \cite{ciocan-fontanine:quantum*1} are
easy consequences of the quantum Monk's formula.  In fact, the quantum
Monk's formula implies that Ciocan-Fontanine's classes satisfy the
same recursive relations as those defining the quantum elementary
symmetric polynomials (cf.~\cite[Lemma~4.2]{postnikov:on*11}).  These
results in turn are the only facts required in the combinatorial proof
of the quantum Giambelli formula for flag varieties given in
\cite{fomin.gelfand.ea:quantum}.  Alternatively, the quantum Schubert
polynomials constructed in \cite{fomin.gelfand.ea:quantum} can easily
be computed by using only the quantum Monk's formula (cf.~\cite[\S
8]{fomin.gelfand.ea:quantum} and \cite[(4.16)]{macdonald:notes}).  The
quantum Pieri formula of Ciocan-Fontanine \cite{ciocan-fontanine:on}
can also be derived combinatorially from the quantum Monk's formula
\cite{postnikov:on*11, fomin.kirillov:quadratic}.  For a survey of
combinatorial approaches to quantum cohomology of flag varieties we
refer the reader to \cite{fomin:lecture}.


In \refsec{sec:schubert} we fix notation regarding Schubert varieties
in partial flag varieties and prove a result which relates the
Schubert varieties in different partial flag varieties.  In
\refsec{sec:curves} we give some tools for handling rational curves in
flag varieties.  The proof of the quantum Monk's formula is finally
given in \refsec{sec:quantum} after a short introduction of the
quantum ring of a flag variety.

We thank S.~Fomin for helpful comments to our paper.


\section{Schubert varieties in partial flag varieties}
\label{sec:schubert}

Set $E = \C^n$ and let $\Fl(E) = \{ V_1 \subset V_2 \subset \dots
\subset V_{n-1} \subset E \mid \dim V_i = i \}$ denote the variety of
full flags in $E$.  Given a fixed flag $F_1 \subset F_2 \subset \dots
\subset F_{n-1} \subset E$ and a permutation $w \in S_n$ there is a
Schubert variety
\[ \Omega_w(F_\bull) = \{ V_\bull \in \Fl(E) \mid \dim(V_p \cap F_q)
   \geq p - r_w(p,n-q) ~\forall p,q \} 
\]
where $r_w(p,q) = \# \{ i \leq p \mid w(i) \leq q \}$.  The
codimension of this variety is equal to the length $\ell(w)$ of $w$.
Notice that the rank conditions on $V_p$ are redundant unless $w$ has
a descent at position $p$, {i.e.\ }$w(p) > w(p+1)$.

Given a sequence of integers $a = (a_1 \leq a_2 \leq \dots \leq a_k)$
with $a_1 \geq 0$ and $a_k \leq n$, we have the partial flag variety
\mbox{$\Fl(a;E) = \{ V_{a_1} \subset \dots \subset V_{a_k} \subset E
  \mid \dim V_{a_i} = a_i \}$.}  Although all such varieties can be
obtained from strictly increasing sequences $a$, it will be convenient
to allow weakly increasing sequences in the notation.  Similarly it is
useful to set $a_0 = 0$ and $a_{k+1} = n$.  Let $S_n(a) \subset S_n$
denote the set of permutations whose descent positions are contained
in the set $\{a_1, a_2, \dots, a_k\}$.  The Schubert varieties in
$\Fl(a;E)$ are indexed by these permutations; the Schubert variety
corresponding to $w \in S_n(a)$ is given by
\[ \Omega_w^{(a)}(F_\bull) = \{ V_\bull \in \Fl(a;E) \mid
   \dim(V_{a_i} \cap F_q) \geq a_i - r_w(a_i,n-q) ~\forall i,q \} \,.
\]
Let $\rho_a : \Fl(E) \to \Fl(a;E)$ be the projection which maps a full
flag $V_\bull$ to the subflag $V_{a_1} \subset \dots \subset V_{a_k}$.
Then for any $w \in S_n(a)$ we have
$\rho_a^{-1}(\Omega^{(a)}_w(F_\bull)) = \Omega_w(F_\bull)$.  On the
other hand, if $w \in S_n$ is any permutation then
$\rho_a(\Omega_w(F_\bull)) = \Omega^{(a)}_{\tilde w}(F_\bull)$ where
$\tilde w \in S_n(a)$ is the permutation obtained from $w$ by
rearranging the elements\linebreak $w(a_i+1), w(a_i+2), \dots, w(a_{i+1})$ in
increasing order for each $0 \leq i \leq k$.  In other words, $\tilde
w$ is the shortest representative for $w$ modulo the subgroup $W_a
\subset S_n$ generated by the simple reflections $s_i = (i,i+1)$ for
$i \not \in \{a_1,\dots,a_k\}$.  For example, if $n = 6$, $a = (2,
5)$, and $w = 6\,2\,3\,1\,5\,4$ then $\tilde w = 2\,6\,1\,3\,5\,4$.

Now let $b = (b_1 \leq b_2 \leq \dots \leq b_k)$ be another sequence
with the same length as $a$, such that $b_i \leq a_i$ for each $i$.
Given a permutation $w \in S_n(a)$ we will need a description of the
set $\{ K_\bull \in \Fl(b;E) \mid \exists~ V_\bull \in
\Omega^{(a)}_w(F_\bull) : K_{b_i} \subset V_{a_i} ~\forall i \}$.

We construct a permutation $\bar w \in S_n(b)$ from $w$ as follows.
Set $w^{(0)} = w$.  Then for each $1 \leq i \leq k$ we let $w^{(i)}$
be the permutation obtained from $w^{(i-1)}$ by rearranging the
elements $w^{(i-1)}(b_i+1), \dots, w^{(i-1)}(a_{i+1})$ in increasing
order.  Finally we set $\bar w = w^{(k)}$.  For example, if $n = 6$,
$a = (2,5)$, $b = (1,2)$, and $w = 2\,6\,3\,4\,5\,1$ then $w^{(1)} =
2\,3\,4\,5\,6\,1$ and $\bar w = 2\,3\,1\,4\,5\,6$.

\begin{lemma} \label{lemma:kernel}
  The set $\{ K_\bull \in \Fl(b;E) \mid \exists~ V_\bull \in
  \Omega^{(a)}_w(F_\bull) : K_{b_i} \subset V_{a_i}
  ~\forall i \}$ is equal to the Schubert variety $\Omega_{\bar
  w}^{(b)}(F_\bull)$ in $\Fl(b;E)$.
\end{lemma}
\begin{proof}
  We prove that the subset $\Omega_i$ of $\Fl_i = \Fl(b_1,\dots,b_i,
  a_{i+1},\dots,a_k; E)$ defined by $\Omega_i = \{ K_\bull \mid
  \exists~ V_\bull \in \Omega^{(a)}_w(F_\bull) : K_{b_j} \subset
  V_{a_j} \text{ for $j \leq i$ and } K_{a_j} = V_{a_j} \text{ for $j
    > i$} \}$ is equal to the Schubert variety in $\Fl_i$ given by the
  permutation $w^{(i)}$.  This is true when $i = 0$.  Let $\rho_j :
  \Fl(E) \to \Fl_j$ denote the projection.  Then it is easy to check
  that $\Omega_{i+1} = \rho_{i+1}(\rho_i^{-1}(\Omega_i))$, so the
  lemma follows from the above remarks about images and inverse images
  of projections $\rho_a$.
\end{proof}

\reflemma{lemma:kernel} has a dual version which we will also need.
Let $a$ and $c$ be weakly increasing sequences of integers between $0$
and $n$, each of length $k$, such that $a_i \leq c_i$ for each $1 \leq
i \leq k$.  Given $w \in S_n(a)$ we define a permutation $\hat w \in
S_n(c)$ as follows.  Set $w^{(k+1)} = w$.  For each $i = k, k-1,
\dots, 1$ we then let $w^{(i)}$ be the permutation obtained from
$w^{(i+1)}$ by rearranging the elements $w^{(i+1)}(a_{i-1}+1), \dots,
w^{(i+1)}(c_i)$ in increasing order.  Finally we set $\hat w = w^{(1)}$.

\begin{lemma} \label{lemma:span}
  The set $\{ W_\bull \in \Fl(c;E) \mid \exists~ V_\bull \in
  \Omega_w^{(a)}(F_\bull) : V_{a_i} \subset W_{c_i}
  ~\forall i \}$ is equal to the Schubert variety $\Omega^{(c)}_{\hat
    w}(F_\bull)$ in $\Fl(c;E)$.
\end{lemma}

Notice that the definitions of the permutations $\bar w$ and $\hat w$
imply that $\ell(\bar w) \geq \ell(w) - \sum_{i=1}^k (a_i -
b_i)(a_{i+1} - a_i)$ and $\ell(\hat w) \geq \ell(w) - \sum_{i=1}^k
(c_i - a_i)(a_i - a_{i-1})$.  In particular, if $a = (1,2,\dots,n-1)$
so that $\Fl(a;E) = \Fl(E)$ then $\ell(\bar w) \geq \ell(w) -
\sum_{i=1}^{n-1} (i - b_i)$ and $\ell(\hat w) \geq \ell(w) -
\sum_{i=1}^{n-1} (c_i - i)$.

\section{Rational curves in partial flag varieties}
\label{sec:curves}

By a rational curve in $\Fl(a;E)$ we will mean the image $C$ of a
regular function $\P^1 \to \Fl(a;E)$.  (We will tolerate that a
rational curve can be a point according to this definition.)  Given a
rational curve $C \subset \Fl(a;E)$ we let \mbox{$C_i = \rho_{a_i}(C)
  \subset \Gr(a_i,E)$} be the image of $C$ by the projection
$\rho_{a_i} : \Fl(a;E) \to \Gr(a_i,E)$.  This curve $C_i$ then has a
{\em kernel\/} and a {\em span\/} \cite{buch:quantum}.  The kernel is
the largest subspace of $E$ contained in all the $a_i$-dimensional
subspaces of $E$ corresponding to points of $C_i$.  We let $b_i$ be
the dimension of this kernel and denote the kernel itself by
$K_{b_i}$.  Similarly, the span of $C_i$ is the smallest subspace of
$E$ containing all the subspaces given by points of $C_i$.  We let
$c_i$ be the dimension of this span and denote the span by $W_{c_i}$.
These subspaces define partial flags $K_\bull \in \Fl(b;E)$ and
$W_\bull \in \Fl(c;E)$ where $b = (b_1\dots,b_k)$ and $c =
(c_1,\dots,c_k)$, which we will call the kernel and span of $C$.

\begin{prop} \label{prop:kernelin}
  Let $C \subset \Fl(a;E)$ be a rational curve with kernel $K_\bull
  \in \Fl(b;E)$ and span $W_\bull \in Fl(c;E)$ and let $w \in S_n(a)$.
  If $C \cap \Omega_w^{(a)}(F_\bull) \neq \emptyset$ then $K_\bull \in
  \Omega_{\bar w}^{(b)}(F_\bull)$ and $W_\bull \in \Omega_{\hat
    w}^{(c)}(F_\bull)$.
\end{prop}
\begin{proof}
  If $V_\bull \in C \cap \Omega_w^{(a)}(F_\bull)$ then we have
  $K_{b_i} \subset V_{a_i} \subset W_{c_i}$ for all $i$.  The
  proposition therefore follows from \reflemma{lemma:kernel} and
  \reflemma{lemma:span}.
\end{proof}

Now let $a = (a_1 < a_2 < \dots < a_k)$ be a strictly increasing
sequence of integers with $1 \leq a_i \leq n-1$.  Recall that the
degree of an algebraic curve $C \subset \P^N$ is the number of points
in the intersection of $C$ with a general hyperplane.  Given a curve
$C \subset \Fl(a;E)$ we let $d_i$ be the degree of $C_i$ in the
Pl{\"u}cker embedding $\Gr(a_i;E) \subset \P(\bigwedge^{a_i}E)$.  The
multidegree of $C$ is defined as the sequence $d = (d_1, \dots, d_k)$.
If $K_\bull \in \Fl(b;E)$ is the kernel and $W_\bull \in \Fl(c;E)$ the
span of $C$, it follows from \linebreak \cite[Lemma 1]{buch:quantum}
that $b_i \geq a_i - d_i$ and $c_i \leq a_i + d_i$ for all $1 \leq i
\leq k$.

Next we shall need a fact about rational curves in the full flag
variety $\Fl(E)$.  For integers $1 \leq i < j \leq n$, let $d_{ij} =
(0,\dots,0, 1,\dots,1, 0,\dots,0)$ denote the multidegree consisting
of $i-1$ zeros followed by $j-i$ ones followed by $n-j$ zeros, {i.e.\ 
  }$(d_{ij})_p = 1$ for $i \leq p < j$ and $(d_{ij})_p = 0$ otherwise.
We set $a = (1,2,\dots,n-1)$ and $b = a - d_{ij} =
(b_1,\dots,b_{n-1})$ where $b_p = p - (d_{ij})_p$.  A partial flag
$K_\bull \in \Fl(b;E)$ together with a subspace $W \subset E$ of
dimension $i+1$ such that $K_{j-2} \cap W = K_{i-1}$ and $K_{j-2} + W
= K_j$ determine a rational curve $C \subset \Fl(E)$ of multidegree
$d_{ij}$, consisting of the flags
\[ V_\bull = ( \, K_1 \subset \dots \subset K_{i-1} \subset L \subset K_i + L
   \subset \dots \subset K_{j-2} + L \subset K_j \subset \dots \subset
   K_{n-1} \, )
\]
for all $i$-dimensional subspaces $L$ such that $K_{i-1} \subset L
\subset W$.  Notice that $K_\bull$ can be recovered as the kernel of
$C$ and $W$ is the subspace of dimension $i+1$ in the span.

\begin{prop} \label{prop:curvedata}
  The curves $C$ in $\Fl(E)$ of multidegree $d_{ij}$ are in 1-1
  correspondence with the pairs $(K_\bull, W)$ where $K_\bull$ is a
  point of $\Fl(b;E)$ and $W \subset E$ is a subspace of dimension
  $i+1$ such that $K_{j-2} \cap W = K_{i-1}$ and $K_{j-2} + W = K_j$.
\end{prop}
\begin{proof}
  Let $C \subset \Fl(E)$ be a rational curve of multidegree $d_{ij}$,
  let $K_\bull \in \Fl(b;E)$ be the kernel of $C$, and let $W_\bull
  \in \Fl(c;E)$ be the span.  It suffices to prove that $W_{i+1} \cap
  K_{j-2} = K_{i-1}$, which in turn is easy to deduce from the
  identities
\[ W_{p+2} \cap K_{p+1} = K_p \]
for all $i-1 \leq p \leq j-3$.  To prove this, let $C$ be the image of
a map $f : \P^1 \to \Fl(E)$, and write $f(x) = (V_1(x) \subset \dots
\subset V_{n-1}(x))$ for $x \in \P^1$.  Since $K_p$ and $W_{p+2}$ are
the kernel and span of the curve $x \mapsto V_{p+1}(x)$ in
$\Gr(p+1,E)$ we see that $K_p \subset W_{p+2}$.  Write $W_{p+2} = K_p
\oplus U$ where $U$ is a 2-dimensional subspace.  Since $K_p \subset
K_{p+1}$ it is enough to show that $U \cap K_{p+1} = 0$.  If this is
false, we can find a basis $\{ u_1, u_2 \}$ for $U$ with $u_1 \in
K_{p+1}$.  Now for suitable coordinates $(s:t)$ on $\P^1$ we have
$V_{p+1}(s:t) = K_p \oplus \C (s u_1 + t u_2)$.  Since $V_{p+2}(s:t)$
contains both $V_{p+1}(s:t)$ and $K_{p+1}$, we conclude that
$V_{p+2}(s:t) = K_{p+1} \oplus \C u_2$ whenever $(s:t) \neq (1:0)$.
In other words the curve $x \mapsto V_{p+2}(x)$ is constant which
contradicts that its degree is one.
\end{proof}

\section{Quantum cohomology of flag varieties}
\label{sec:quantum}

For each permutation $w \in S_n$ we let $\Omega_w$ denote the class of
$\Omega_w(F_\bull)$ in the cohomology ring $H^* \Fl(E) = H^*(\Fl(E);
\Z)$.  The Schubert classes $\Omega_w$ form a basis for this ring.  If
$d = (d_1,\dots,d_{n-1})$ is a multidegree we set $|d| = \sum d_i$.
Given three permutations $u, v, w \in S_n$ such that $\ell(u) +
\ell(v) + \ell(w) = \binom{n}{2} + 2 |d|$, the Gromov-Witten invariant
$\left< \Omega_u, \Omega_v, \Omega_w \right>_d$ is defined as the
number of rational curves of multidegree $d$ in $\Fl(E)$ which meet
each of the Schubert varieties $\Omega_u(F_\bull)$,
$\Omega_v(G_\bull)$, and $\Omega_w(H_\bull)$ for general fixed flags
$F_\bull, G_\bull, H_\bull$ in $E$.  If $\ell(u) + \ell(v) + \ell(w)
\neq \binom{n}{2} + 2 |d|$ then $\left< \Omega_u, \Omega_v, \Omega_w
\right>_d = 0$.

Let $q_1, \dots, q_{n-1}$ be independent variables, and write
$\Z[q] = \Z[q_1,\dots,q_{n-1}]$.  The quantum cohomology ring $\QH^*
\Fl(E)$ is a $\Z[q]$-algebra which is isomorphic to $H^* \Fl(E)
\otimes \Z[q]$ as a module over $\Z[q]$.  In this ring we have quantum
Schubert classes $\sigma_w = \Omega_w \otimes 1$.  Multiplication in
$\QH^* \Fl(E)$ is defined by
\[ \sigma_u \cdot \sigma_v = \sum_{w,d} \left< \Omega_u, \Omega_v,
  \Omega_{w^\vee} \right>_d \, q^d \,\sigma_w 
\]
where the sum is over all permutations $w \in S_n$ and multidegrees $d
= (d_1,\dots,d_{n-1})$; here we set $q^d = \prod q_i^{d_i}$ and we let
$w^\vee \in S_n$ denote the permutation of the dual Schubert class to
$\Omega_w$, {i.e.\ }$w^\vee = w_0 w$ where $w_0$ is the longest
permutation in $S_n$.  It is a non-trivial fact that this defines an
associative ring \cite{ruan.tian:mathematical,
  kontsevich.manin:gromov-witten, fulton.pandharipande:notes}.

For $1 \leq i < j \leq n$ we let $t_{ij} = (i,j) \in S_n$ denote the
transposition which interchanges $i$ and $j$.  We furthermore set
$q_{ij} = q^{d_{ij}} = q_i\, q_{i+1} \dots q_{j-1}$.  Our goal is to
prove the following quantum version of the Monk's formula from
\cite{fomin.gelfand.ea:quantum}.

\begin{thm}
For $w \in S_n$ and $1 \leq r < n$ we have
\[ \sigma_{s_r} \cdot \sigma_w = \sum \sigma_{w\, t_{kl}} + 
   \sum q_{ij}\, \sigma_{w\, t_{ij}}
\]
where the first sum is over all transpositions $t_{kl}$ such that
$k \leq r < l$ and $\ell(w\, t_{kl}) = \ell(w) + 1$, and the second
sum is over all transpositions $t_{ij}$ such that $i \leq r < j$
and $\ell(w\, t_{ij}) = \ell(w) - \ell(t_{ij}) = \ell(w) -
2(j-i) + 1$.
\end{thm}
\begin{proof}
  The first sum is dictated by the classical Monk's formula
  \cite{monk:geometry}.  The second sum is equivalent to the following
  statement.  If $d = (d_1,\dots,d_{n-1})$ is a non-zero multidegree
  and $u, w \in S_n$ are permutations such that $\ell(u) + \ell(w) +
  \ell(s_r) = \binom{n}{2} + 2 |d|$ then the Gromov-Witten invariant
  $\left< \Omega_u, \Omega_w, \Omega_{s_r} \right>_d$ is equal to one
  if $d = d_{ij}$ for some $i,j$ such that $i \leq r < j$ and $u^{-1}
  w_0 w = t_{ij}$; otherwise $\left< \Omega_u, \Omega_w, \Omega_{s_r}
  \right>_d = 0$.
  
  Suppose $\left< \Omega_u, \Omega_w, \Omega_{s_r} \right>_d \neq 0$
  and let $C$ be a rational curve of multidegree $d$ which meets three
  Schubert varieties $\Omega_u(F_\bull)$, $\Omega_w(G_\bull)$, and
  $\Omega_{s_r}(H_\bull)$ in general position.  Let $K_\bull \in
  \Fl(b;E)$ be the kernel of $C$ and set $a = (1,2,\dots,n-1)$.  Then
  $b_p \geq a_p - d_p$ for all $1 \leq p \leq n-1$.  By
  \refprop{prop:kernelin} we have $K_\bull \in \Omega_{\bar
    u}^{(b)}(F_\bull) \cap \Omega_{\bar w}^{(b)}(G_\bull) \cap
  \Omega_{\bar s_r}^{(b)}(H_\bull)$.  Since the flags are general this
  implies that $\ell(\bar u) + \ell(\bar w) + \ell(\bar s_r) \leq \dim
  \Fl(b;E)$.  On the other hand the inequalities $\ell(\bar u) \geq
  \ell(u) - \sum (p - b_p)$, $\ell(\bar w) \geq \ell(w) - \sum (p -
  b_p)$, $\ell(\bar s_r) \geq 0$, and $\sum (p - b_p) \leq |d|$ imply
  that $\ell(\bar u) + \ell(\bar w) + \ell(\bar s_r) \geq \binom{n}{2}
  - 1$.  Since this is the maximal possible dimension of $\Fl(b;E)$ we
  conclude that all inequalities are satisfied with equality.
  
  This first implies that $b = a-d = (1-d_1, 2-d_2, \dots,
  n-1-d_{n-1})$.  Furthermore, since $\dim \Fl(b;E) = \binom{n}{2} -
  1$ we deduce that $d = d_{ij}$ for some $1 \leq i < j \leq n$.  Thus
  $\Fl(b;E) = \Fl(1,\dots, j-2, j, \dots,n-1; E)$ is the variety of
  partial flags with subspaces of all dimensions other than $j-1$.
  Since $\ell(\bar s_r) = 0$ it follows that $i \leq r < j$.  The fact
  that $\ell(\bar u) = \ell(u) - |d|$ implies that $\bar u = u \, s_i
  s_{i+1} \cdots s_{j-1}$ by the definition of $\bar u$.  Similarly we
  have $\bar w = w \, s_i s_{i+1} \cdots s_{j-1}$.  Now since
  $\ell(\bar u) + \ell(\bar w) = \dim \Fl(b;E)$ and $\Omega_{\bar
    u}^{(b)}(F_\bull) \cap \Omega_{\bar w}^{(b)}(G_\bull) \neq
  \emptyset$ we conclude that $\bar u$ and $\bar w$ are dual with
  respect to $\Fl(b;E)$, {i.e.\ }$\bar u^{-1} w_0 \bar w = s_{j-1}$ or
  equivalently $u^{-1} w_0 w = t_{ij}$ as required.
  
  It remains to be proven that if $d = d_{ij}$ and $u^{-1} w_0 w =
  t_{ij}$ for some $i \leq r < j$ then there exists a unique rational
  curve of multidegree $d$ which meets the three given Schubert
  varieties.  Set $\bar u = u \, s_i s_{i+1} \cdots s_{j-1}$ and $\bar
  w = w \, s_i s_{i+1} \cdots s_{j-1}$.  Since $\ell(u \, t_{ij}) =
  \ell(w_0 w) = \binom{n}{2} - \ell(w) = \ell(u) - \ell(t_{ij})$ it
  follows that $\ell(\bar u) = \ell(u) - |d|$ and similarly $\ell(\bar
  w) = \ell(w) - |d|$.  Thus $\ell(\bar u) + \ell(\bar w) = \dim
  \Fl(b;E)$ where $b = a - d$.  Since $\bar u^{-1} w_0 \bar w =
  s_{j-1}$ we conclude that there is a unique partial flag $K_\bull
  \in \Omega_{\bar u}^{(b)}(F_\bull) \cap \Omega_{\bar
    w}^{(b)}(G_\bull)$.  Similarly, if we set $\hat u = u \, s_{j-1}
  s_{j-2} \cdots s_i$ and $\hat w = w \, s_{j-1} s_{j-2} \cdots s_i$
  then there exists a unique partial flag $W_\bull \in \Omega_{\hat
    u}^{(c)}(F_\bull) \cap \Omega_{\hat w}^{(c)}(G_\bull)$ where $c =
  a+d$.
  
  In fact, we can say precisely what these partial flags look like.
  For $1 \leq p \leq n$ we set $L_p = F_{n+1-p} \cap G_p$.  Since the
  flags $F_\bull$ and $G_\bull$ are general, these spaces have
  dimension one, and $E = L_1 \oplus \dots \oplus L_n$.  Now $K_p =
  L_{\bar u(1)} \oplus \cdots \oplus L_{\bar u(p)}$ for each $p \neq
  j-1$ and $W_p = L_{\hat u(1)} \oplus \cdots \oplus L_{\hat u(p)}$
  for $p \neq i$.  Otherwise stated we have $K_p = W_p = L_{u(1)}
  \oplus \cdots \oplus L_{u(p)}$ for $1 \leq p \leq i-1$ and for $j
  \leq p < n$.  For $i-1 \leq p \leq j-2$ we have $K_p = K_{i-1}
  \oplus L_{u(i+1)} \oplus \cdots \oplus L_{u(p+1)}$ while $W_{p+2} =
  K_p \oplus U$ where $U = L_{u(i)} \oplus L_{u(j)}$.  In particular
  we get $W_{i+1} \cap K_{j-2} = K_{i-1}$ and $W_{i+1} + W_{j-2} =
  W_j$ so by \refprop{prop:curvedata} there is exactly one rational
  curve of multidegree $d$ with kernel $K_\bull$ and span $W_\bull$.
  This curve consists of all flags
\[ V_\bull = ( \, K_1 \subset \dots \subset K_{i-1} \subset K_{i-1} \oplus
  L \subset \dots \subset K_{j-2} \oplus L \subset K_j \subset \dots
  \subset K_{n-1} \, )
\]
where $L \subset U$ is a one dimensional subspace.  When $L =
L_{u(i)}$ we have $V_\bull \in \Omega_{u}(F_\bull)$, while $V_\bull
\in \Omega_w(G_\bull)$ when $L = L_{u(j)}$.  Finally, $V_\bull$
belongs to $\Omega_{s_r}(H_\bull)$ if and only if $V_r \cap H_{n-r}
\neq 0$.  Now take any non-zero element $x \in W_{r+1} \cap H_{n-r}$
and let $x'$ be the $U$-component of $x$ in $W_{r+1} = K_{r-1} \oplus
U$.  Taking $L = \C \, x'$ then gives a point $V_\bull \in
\Omega_{s_r}(H_\bull)$.  This completes the proof.
\end{proof}

\newpage


\providecommand{\bysame}{\leavevmode\hbox to3em{\hrulefill}\thinspace}


\end{document}